\documentclass[A4j,11pt]{article}
\usepackage{amssymb,amsmath,amsthm,amscd,amsfonts}
\usepackage{graphicx}
\usepackage{graphics}

\begin{document}

\title{{\bf The preservability \\
of the curvature-adaptedness\\
along the mean curvature flow}}
\author{Naoyuki Koike}
\date{}

\maketitle

%\address{Department of Mathematics, Faculty of Science, Tokyo University of Science, 1-3 Kagurazaka, Shinjuku-ku, Tokyo 162-8601, Japan}
%\email{koike@rs.tus.ac.jp}

\begin{abstract}
In this paper, we investigate the preservability of the curvature-adaptedness along the mean curvature flow 
starting from a compact curvature-adapted hypersurface in locally symmetric spaces, where the curvature-adaptedness means that 
the shape operator and the normal Jacobi operator of the hypersurface commute.  
\end{abstract}

%\keywords{mean curvature flow, curvature-adaptedness, shape operator, normal Jacobi operator}
%\subjclass[2010]{53C40, 53C44}

\section{Introduction}
In 1996, K. Smoczyk (\cite{S}) proved that the Lagrangianity is preserved along the mean curvature flow starting 
from a compact Lagrangian submanifold in a Calabi-Yau manifold.  
He proved this fact by showing that the squared norms of the $2$-forms on the submanifolds induced from 
the fundamental $2$-form of the Calabi-Yau manifold remain to vanish.  
%From this result, the mean curvature flow starting from a Lagrangian submanifold in a Calabi-Yau manifold 
%has been called the {\it Lagrangian mean curvature flow}.  
%See [Sm-1,2,3], [Wang-123] and [Thomas-Yau] etc. about the study of the Lagrangian mean curvature flow.  
In this paper, we investigate the preservability of the curvature-adaptedness along the mean curvature flow 
starting from a compact curvature-adapted hypersurface in locally symmetric spaces, where the curvature-adaptedness of 
a hypersurface means that the shape operator and the normal Jacobi operator of the hypersurface commute.  
See the following paragraph about the precise definition of the curvature-adaptedness.  
The strategy of the proof of the main theorem (Theorem B) of this paper is to find a sufficient condition for the squared norm of 
the commutator of the shape operator and the normal Jacobi operator to remains to vanish by investigating the evolution of the squared norm.  
Our strategy is similar to that of \cite{S}.  

Throughout this paper, we assume that all manifolds are oriented.  
We shall state the definition of the curvature-adaptedness of the hypersurface.  
Let $(\widetilde M,\widetilde g)$ be an $(n+1)$-dimensional Riemannian manifold, $M$ be a $n$-dimensional compact manifold 
and $f$ be an immersion of $M$ into $\widetilde M$.  
Also, let $\xi$ be the unit normal vector field of $f$ compatible with the orientations of $M$ and 
$\widetilde M$.  Denote by $\widetilde R$ the curvature tensor of $\widetilde M$.  
Also, denote by $T^{\perp}_xM$ the normal space of $f$ at $x(\in M)$, by $A$ the shape tensor 
of $f$ for $-\xi$ and by $\nabla^{\perp}$ the normal connection of $f$.  
If the shape operator $A$ commutes with the normal Jacobi operator 
$\widetilde R(\xi):=f_{\ast}^{-1}\circ \widetilde R(\cdot,\xi)\xi\circ f_{\ast}$, then $f$ is said to be {\it curvature-adapted}.  
This notion was introduced by J. Berndt and L. Vanhecke (\cite{BV}).  
All hypersurfaces in real space forms are curvature-adapted and the curvature-adapted hypersurfaces in 
the a complex projective space and the complex hyperbolic space are called Hopf hypersurfaces.  
Note that the condition of the curvature-adapetdness is stricter as the curvature of the ambient space is more complicate.  

Next we shall state the definition of an isoparametric hypersurface.  
A ($C^{\infty}$-)function $\psi$ over $\widetilde M$ is said to be {\it isoparametric} if it satisfied 
the following conditions:

\vspace{0.3truecm}

(i) $\vert\vert d\psi\vert\vert^2=a(\psi)$ holds for some $C^{\infty}$-function $a$ over ${\Bbb R}$, 

(ii) $\widetilde{\triangle}\psi=b(\psi)$ holds for some continuous function $b$ over ${\Bbb R}$, 

where $\widetilde{\triangle}$ is the Laplace-Beltrami operator of $\widetilde M$.  

\vspace{0.3truecm}

\noindent
The regular level sets of an isoprarametric function are called {\it isoparametric hypersurfaces} in 
$\widetilde M$.  Since the regular level sets of an isoprarametric function are parallel to one another and 
they are of constant mean curvature, the mean curvature flow starting from an isoparametric hypersurface 
consists of the parallel hypersurfaces of the hypersurface.  
In this paper, if a hypersurface $M$ and the parallel hypersurfaces sufficiently close to $M$ are of constant mean curvature, then 
we call $M$ a {\it locally isoparametric hypersurface}.  

In 2012, T. Murphy (\cite{M}) studied curvature-adapted hypersurfaces in a compact symmetric space.  
He proved that a curvature-adapted hypersurface in a compact symmetric space is isoparametric if and only if 
both the shape operator and the normal Jacobi operator of the hypersurface have constant eigenvalues.  
In 2014, the author (\cite{K2}) proved that curvature-adapted submanifolds with maximal flat section 
in a symmetric space are principal orbits of the isotropy action of the symmetric space under certain conditions, where 
``submanifold with maximal flat section'' means that the normal umbrellas of the submanifold are maximal dimensional flat totally geodesic 
submanifolds in the ambient symmetric space.  
Since principal orbits of the isotropy action are curvature-adapted isoparametric submanifolds and they are parallel to one another, 
the mean curvature flow starting from a principal orbit of the isotropy action consists of principal orbits of the isotropy action.  
Hence the curvature-adaptedness is preserved along the mean curvature flow starting from the principal orbit.  
From these facts, we can derive that the curvature-adaptedness is preserved along the mean curvature flow starting from 
a curvature-adapted submanifold with maximal flat section satisfying the conditions.  
%%See [AB] about the study (based on the theory of connections on vector bundles) of the mean curvature flow 
%%starting from a submanifold of higher codimension in a general framework.  

Next we shall state the definitions of the mean curvature flow and the backward mean curvature flow.  
Let $M,(\widetilde M,\widetilde g)$ and $f$ be as above, and $\{f_t\}_{t\in[0,T)}$ be a $C^{\infty}$-family 
of immersions of $M$ into $\widetilde M$ and $\xi_t$ be the unit normal vector field of $f_t$ compatible 
with the orietations of $M$ and $\widetilde M$.  
Denote by $H_t$ the mean curvature of $f_t$ for $-\xi_t$.  
Define a map $F:M\times[0,T)\to V$ by $F(x,t):=f_t(x)$ ($(x,t)\in M\times[0,T)$).  
This family $\{f_t\}_{t\in[0,T)}$ is called the {\it mean curvature flow starting from} $f$ 
if $f_0=f$ and if the following evolution equation hold:
$$\frac{\partial F}{\partial t}=-H_t\xi_t.\leqno{(1.1)}$$
Also, this family $\{f_t\}_{t\in[0,T)}$ is called the {\it backward mean curvature flow starting from} $f$ 
if $f_0=f$ and if the following evolution equation hold:
$$\frac{\partial F}{\partial t}=H_t\xi_t.\leqno{(1.2)}$$
Note that G. Huisken (\cite{Hu1,Hu2}) initiated the study of the mean curvature flow starting from a hypersurface as 
the evolution of immersions.  Also, B. Andrews and C. Baker (\cite{AB}) studied the mean curvature flow starting from 
a submanifold (of general codimension) as the evolution of immersions in the aspect of the theory of the vector bundle.  

In this paper, we shall tackle the following question:

\vspace{0.3truecm}

\noindent
{\bf Question 1.} {\sl In what case, is the curvature-adaptedness preserved along the mean curvature flow starting from $f$ 
if $f$ is curvature-adapted?}

\vspace{0.3truecm}

First, we derive the following result for this question.  

\vspace{0.5truecm}

\noindent
{\bf Theorem A.} {\sl 
Assume that $(\widetilde M,\widetilde g)$ is an $(n+1)$-dimensional locally symmetric space.  
Let $f$ be an immersion of an $n$-dimensional compact manifold $M$ into $\widetilde M$ and $\{f_t\}_{t\in[0,T)}$ the mean curvature flow 
starting from $f$.  If $f$ is curvature-adapted and if both the shape operator and the normal Jacobi operator of $f$ have constant eigenvalues, 
then $f$ is locally isoprarametric and $f_t$ remains to be curvature-adapted and locally isoparametric for all $t\in[0,T)$.}

\vspace{0.5truecm}

We shall prepare some notations to state the main result in this paper.  
Let $f$ be an immersion of an $n$-dimensional compact manifold $M$ into a $(n+1)$-dimensional Riemannian manifold $\widetilde M$ and 
$\{f_t\}_{t\in[0,T)}$ the mean curvature flow starting from $f$.  Denote by $g_t$ the induced metric on $M$ by $f_t$, and by $\nabla^t$ and $R_t$ the Levi-Civita connection 
and the curvature tensor of $g_t$, respectively.  
Set $S_t:=[A_t,\widetilde R(\xi_t)]$ (the commutator of $A_t$ and $\widetilde R(\xi_t)$) 
and define non-negative functions $\rho_t$ 
($0\leq t<T$) over $M$ by 
$$\rho_t:=-{\rm Tr}(S_t^2)(={\rm Tr}(S_t^{\ast_t}\circ S_t)=\vert\vert S_t\vert\vert_t^2(\geq 0)),
\leqno{(1.3)}$$
where $A_t$ and $\widetilde R(\xi_t)$ are the shape operator and the normal Jacobi operator of $f_t$, 
respectively, and $S_t^{\ast_t}$ is the adjoint operator of $S_t$ with respect to $g_t$.  
This function $\rho_t$ implies the gap from the curvature-adaptedness of $f_t$.  
In this paper, we call this function a {\it gap function}.  
Define $(1,2)$-tensor fields $\widetilde R_i(\xi_t)$ ($i=1,3$) over $M$ by 
$$\widetilde R_1(\xi_t)(X,Y):=f_{t\ast}^{-1}((\widetilde R(\xi_t,f_{t\ast}X)f_{t\ast}Y)_T)\qquad\,\,
(X,Y\in TM)$$
and 
$$\widetilde R_3(\xi_t)(X,Y):=f_{t\ast}^{-1}(\widetilde R(f_{t\ast}X,f_{t\ast}Y)\xi_t)\qquad\,\,(X,Y\in TM),$$
where $(\bullet)_T$ is the $f_{t\ast}(TM)$-component of $(\bullet)$.  
Also, for tangent vector fields $X$ and $Y$ on $M$, define a $(1,1)$-tensor field $R_t(X,Y)A_t^k$ ($k=1,2$) 
on $M$ by 
$$R_t(X,Y)A_t^k:=\nabla^t_X(\nabla^t_Y(A_t^k))-\nabla^t_Y(\nabla^t_X(A_t^k))-\nabla^t_{[X,Y]}(A_t^k).$$

Define skew-symmetric $(1,1)$-tensor fields $\widehat S_t$ ($0\leq t<T$) over $M$ by 
$$\begin{array}{l}
\displaystyle{\widehat S_t:=2\left[A_t^2+\widetilde R(\xi_t),
{\rm Tr}_{g_t}^{\bullet}(R_t(\cdot,\bullet)A_t)(\bullet))\right]}\\
\hspace{1truecm}\displaystyle{+2\left[A_t,\,(\widetilde R_3(\xi_t)-\widetilde R_1(\xi_t))
(\cdot,{\rm Tr}_{g_t}^{\bullet}(\nabla^t_{\bullet}A_t)(\bullet))\right]}\\
%\hspace{1truecm}\displaystyle{+\left[A_t,\,(\widetilde R_3(\xi_t)-\widetilde R_1(\xi_t))
%(\cdot,\widetilde{\rm Ric}^{\sharp}(\xi_t)_T)\right]}\\
%\hspace{1truecm}\displaystyle{+\left[A_t,\,(\widetilde R_3(\xi_t)-\widetilde R_1(\xi_t))
%(\cdot,{\rm grad}_{g_t}H_t+{\rm Tr}^{\bullet}_{g_t}(\nabla^t_{\bullet}A_t)(\bullet))\right]}\\
\hspace{1truecm}\displaystyle{+2{\rm Tr}^{\bullet}_{g_t}\left[\nabla^t_{\bullet}A_t,\,
\nabla^t_{\bullet}\widetilde R(\xi_t)\right].}
\end{array}\leqno{(1.4)}$$
%where $\widetilde{\rm Ric}^{\sharp}(\xi_t)_T$ denotes the inverse image of the 
%$(f_t)_{\ast}(TM)$-component of $\widetilde{\rm Ric}^{\sharp}(\xi_t)$ by $(f_t)_{\ast}$.  
Note that, if ${\rm Ker}\,A_t=\{0\}$, then the tensor field $\widetilde R_3(\xi_t)-\widetilde R_1(\xi_t)$ 
is described in terms of $\nabla^t\widetilde R(\xi_t)$ and $A_t^{-1}$ (see $(3.8)$) and hence $\widehat S_t$ 
is described in terms of $A_t,A_t^{-1},\nabla^tA_t,\nabla^t\nabla^tA_t,\widetilde R(\xi_t)$ and 
$\nabla^t\widetilde R(\xi_t)$.  
Define a function $\mu_t:M\to\mathbb R$ ($t\in[0,T)$) by 
$$\mu_t(x):=\left\{
\begin{array}{ll}
\displaystyle{-\frac{\langle(\widehat S_t)_x,(S_t)_x\rangle}{\|(S_t)_x\|^2}} & (\|(S_t)_x\|\not=0)\\
\displaystyle{0} & (\|(S_t)_x\|=0),
\end{array}\right.\leqno{(1.5)}$$
where $\langle(\widehat S_t)_x,(S_t)_x\rangle$ denotes 
${\rm Tr}((\widehat S_t^{\ast_t})_x\circ(S_t)_x)(=-{\rm Tr}((\widehat S_t)_x\circ(S_t)_x))$.  

\vspace{0.5truecm}

We prove the following result for Question 1.  
%%preservability theorem of the curvature-adaptedness along the mean curvature flow.  

\vspace{0.5truecm}

\noindent
{\bf Theorem B.} {\sl Let $(\widetilde M,\widetilde g)$ be an $(n+1)$-dimensional locally symmetric space, 
$f$ a curvature-adapted immersion of an $n$-dimensional compact manifold $M$ into $\widetilde M$, 
$\widehat S$ the skew-symmetric $(1,1)$-tensor field on $M$ defined as in $(1.4)$ for $f$ and $\{f_t\}_{t\in[0,T)}$ 
the mean curvature flow starting from $f$.  If $\widehat S\not=0$, then $f_t$ ($0<t<\varepsilon$) are not 
curvature-adapted for some $\varepsilon>0$.}

\vspace{0.5truecm}

From Theorems A and B, we obtain the following result for this question.  

\vspace{0.5truecm}

\noindent
{\bf Corollay C.} {\sl Let $(\widetilde M,\widetilde g)$ be an $(n+1)$-dimensional locally symmetric space, 
$f$ an immersion of an $n$-dimensional compact manifold $M$ into $\widetilde M$ and $\widehat S$ the skew-symmetric $(1,1)$-tensor filed 
on $M$ defined as in $(1.4)$ for $f$.  Then, if $f$ is curvature-adapted and if both the shape operator and the normal Jacobi operator of $f$ 
have constant eigenvalues, then $\widehat S=0$ holds.}

\vspace{0.5truecm}

Naturally the following question arises.  

\vspace{0.5truecm}

\noindent
{\bf Question 2.} {\sl Let $f$ and $\{f_t\}_{t\in[0,T)}$ be as in Theorem B.  Does $f_t$ remain to be curvature-adapted 
for all $t\in[0,T)$ if $\widehat S=0$?}

\vspace{0.5truecm}

If this question were solved affirmatively, then we see that $\widehat S$ is an obstruction for the curvature-adaptedness to be preserved 
along the mean curvature flow starting from a curvature-adapted compact hypersurface in a locally symmetric space.  

We prove the following result for this question.  
%%preservability theorem of the curvature-adaptedness along the mean curvature flow.  

\vspace{0.5truecm}

\noindent
{\bf Theorem D.} {\sl Let $(\widetilde M,\widetilde g)$ be an $(n+1)$-dimensional locally symmetric space, 
$f$ a curvatue-adapted immersion of an $n$-dimensional compact manifold $M$ into $\widetilde M$ and $\{f_t\}_{t\in[0,T)}$ 
the mean curvature flow starting from $f$.  If $\widehat S_0=0$ and if 
$$\mathop{\sup}_{t\in[0,T)}\,\mathop{\sup}_{x\in M}\,\mu_t(x)<\infty,$$
then $f_t$ remains to be curvature-adapted for all $t\in[0,T)$.}

\vspace{0.5truecm}

As a corollary of Theorem D, we obtain the following result.  

\vspace{0.5truecm}

\noindent
{\bf Corollary E.} {\sl Let $(\widetilde M,\widetilde g)$ be an $(n+1)$-dimensional locally symmetric space, 
$f$ a curvatue-adapted immersion of an $n$-dimensional compact manifold $M$ into $\widetilde M$ and $\{f_t\}_{t\in[0,T)}$ 
the mean curvature flow starting from $f$.  If $\widehat S_0=0$ and if 
$$\mathop{\inf}_{t\in[0,T)}\,\mathop{\min}_{x\in M}\,\langle(\widehat S_t)_x,(S_t)_x\rangle\geq 0,$$
then $f_t$ remains to be curvature-adapted for all $t\in[0,T)$.}

\vspace{0.5truecm}

Assume that $f$ is curvature-adapted and that $f_{t_0}$ is not curvature-adapted for some $t_0\in[0,T)$.  Set 
$$t_{\min}:=\inf\,\{t\in[0,T)\,|\,f_t\,\,{\rm is}\,\,{\rm not}\,\,{\rm curvature-adapted}\}.$$
Then, according to Theorem D, $\displaystyle{\mathop{\sup}_{x\in M}\,\mu_t(x)}$ diverges to $+\infty$ 
as $t\downarrow t_{\min}$ (see Figure 1.1).  This fact is restated in terms of ``backward mean curvature flow'' as follows.  

\vspace{0.5truecm}

\noindent
{\bf Theorem F.} {\sl Let $(\widetilde M,\widetilde g)$ be an $(n+1)$-dimensional locally symmetric space, 
$f$ an immersion of an $n$-dimensional compact manifold $M$ into $\widetilde M$ and $\{f^b_t\}_{t\in[0,T)}$ 
the backward mean curvature flow starting from $f$.  Assume that $f$ is not curvature-adapted and that 
$f^b_{t_0}$ is curvature-adapted for some $t_0\in[0,T)$, where $t_0$ is the first time such that $f_t^b$ is curvature-adapted.  
%%$t_0:=\min\{t\in[0,T)\,|\,f_t\,\,{\rm is}\,\,{\rm curvature-adapted}\}$.  
Then $\displaystyle{\lim_{t\uparrow t_0}\,\mathop{\sup}_{x\in M}\,\mu_t(x)=\infty}$ holds.}

\vspace{0.5truecm}

\centerline{
%WinTpicVersion3.08
\unitlength 0.1in
\begin{picture}( 73.0000, 21.7000)(-32.3000,-27.1000)
% STR 2 0 3 0
% 3 2620 2610 2620 2710 1 0
% $\,$
\put(26.2000,-27.1000){\makebox(0,0)[lt]{$\,$}}%
% VECTOR 2 0 3 0
% 2 1600 2400 4000 2400
% 
\special{pn 8}%
\special{pa 1600 2400}%
\special{pa 4000 2400}%
\special{fp}%
\special{sh 1}%
\special{pa 4000 2400}%
\special{pa 3934 2380}%
\special{pa 3948 2400}%
\special{pa 3934 2420}%
\special{pa 4000 2400}%
\special{fp}%
% VECTOR 2 0 3 0
% 2 1800 2600 1800 600
% 
\special{pn 8}%
\special{pa 1800 2600}%
\special{pa 1800 600}%
\special{fp}%
\special{sh 1}%
\special{pa 1800 600}%
\special{pa 1780 668}%
\special{pa 1800 654}%
\special{pa 1820 668}%
\special{pa 1800 600}%
\special{fp}%
% DOT 0 0 3 0
% 2 2280 2400 2280 2400
% 
\special{pn 20}%
\special{sh 1}%
\special{ar 2280 2400 10 10 0  6.28318530717959E+0000}%
\special{sh 1}%
\special{ar 2280 2400 10 10 0  6.28318530717959E+0000}%
% LINE 1 0 3 0
% 2 1800 2400 2270 2400
% 
\special{pn 13}%
\special{pa 1800 2400}%
\special{pa 2270 2400}%
\special{fp}%
% ELLIPSE 1 0 3 0
% 4 2280 390 4480 2400 2280 3380 4120 2400
% 
\special{pn 13}%
\special{ar 2280 390 2200 2010  0.8742723 1.5707963}%
% LINE 2 2 3 0
% 2 3700 1930 3700 2400
% 
\special{pn 8}%
\special{pa 3700 1930}%
\special{pa 3700 2400}%
\special{dt 0.045}%
% DOT 0 0 3 0
% 2 3700 2400 3700 2400
% 
\special{pn 20}%
\special{sh 1}%
\special{ar 3700 2400 10 10 0  6.28318530717959E+0000}%
\special{sh 1}%
\special{ar 3700 2400 10 10 0  6.28318530717959E+0000}%
% LINE 2 2 3 0
% 2 2280 2390 2280 620
% 
\special{pn 8}%
\special{pa 2280 2390}%
\special{pa 2280 620}%
\special{dt 0.045}%
% LINE 2 2 3 0
% 2 3700 1930 3700 1670
% 
\special{pn 8}%
\special{pa 3700 1930}%
\special{pa 3700 1670}%
\special{dt 0.045}%
% STR 2 0 3 0
% 3 4070 2260 4070 2360 1 0
% $t$
\put(40.7000,-23.6000){\makebox(0,0)[lt]{$t$}}%
% STR 2 0 3 0
% 3 2220 2360 2220 2460 1 0
% $t_{\min}$
\put(22.2000,-24.6000){\makebox(0,0)[lt]{$t_{\min}$}}%
% STR 2 0 3 0
% 3 3640 2360 3640 2460 1 0
% $T$
\put(36.4000,-24.6000){\makebox(0,0)[lt]{$T$}}%
% STR 2 0 3 0
% 3 1710 610 1710 710 3 0
% $y$
\put(17.1000,-7.1000){\makebox(0,0)[rb]{$y$}}%
% VECTOR 2 1 3 0
% 2 1510 1720 2030 2390
% 
\special{pn 8}%
\special{pa 1510 1720}%
\special{pa 2030 2390}%
\special{da 0.070}%
\special{sh 1}%
\special{pa 2030 2390}%
\special{pa 2006 2326}%
\special{pa 1998 2348}%
\special{pa 1974 2350}%
\special{pa 2030 2390}%
\special{fp}%
% VECTOR 2 1 3 0
% 2 1510 1720 3060 2260
% 
\special{pn 8}%
\special{pa 1510 1720}%
\special{pa 3060 2260}%
\special{da 0.070}%
\special{sh 1}%
\special{pa 3060 2260}%
\special{pa 3004 2220}%
\special{pa 3010 2242}%
\special{pa 2990 2258}%
\special{pa 3060 2260}%
\special{fp}%
% VECTOR 2 1 3 0
% 2 1500 1220 2130 2390
% 
\special{pn 8}%
\special{pa 1500 1220}%
\special{pa 2130 2390}%
\special{da 0.070}%
\special{sh 1}%
\special{pa 2130 2390}%
\special{pa 2116 2322}%
\special{pa 2106 2344}%
\special{pa 2082 2342}%
\special{pa 2130 2390}%
\special{fp}%
% VECTOR 2 1 3 0
% 2 1500 1220 2680 1460
% 
\special{pn 8}%
\special{pa 1500 1220}%
\special{pa 2680 1460}%
\special{da 0.070}%
\special{sh 1}%
\special{pa 2680 1460}%
\special{pa 2620 1428}%
\special{pa 2628 1450}%
\special{pa 2612 1466}%
\special{pa 2680 1460}%
\special{fp}%
% STR 2 0 3 0
% 3 1450 1670 1450 1770 3 0
% $\displaystyle{y=\mathop{\max}_{x\in M}\,\rho_t(x)}$
\put(14.5000,-17.7000){\makebox(0,0)[rb]{$\displaystyle{y=\mathop{\max}_{x\in M}\,\rho_t(x)}$}}%
% STR 2 0 3 0
% 3 1450 1160 1450 1260 3 0
% $\displaystyle{y=\mathop{\sup}_{x\in M}\,\mu_t(x)}$
\put(14.5000,-12.6000){\makebox(0,0)[rb]{$\displaystyle{y=\mathop{\sup}_{x\in M}\,\mu_t(x)}$}}%
% ELLIPSE 1 0 3 0
% 4 3360 570 4400 1740 2010 650 3770 1870
% 
\special{pn 13}%
\special{ar 3360 570 1040 1170  1.2299661 3.0890485}%
\end{picture}%
\hspace{9.5truecm}}

\vspace{0.5truecm}

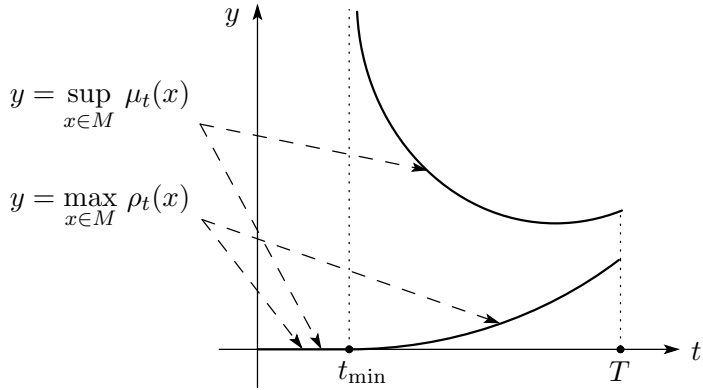
\centerline{{\bf Figure 1.1: {\small The graph of $\displaystyle{\mathop{\sup}_{x\in M}\,\mu_t(x)}$}}}

\vspace{0.5truecm}

In the future, we plan to tackle the following question.  

\vspace{0.5truecm}

\noindent
{\bf Question 3.}\ {\sl Can we find a pinching condition of the norms $\|A_0\|,\,\|\nabla^0A_0\|$ and $\|(\nabla^0)^2A_0\|$ 
satisfying $\displaystyle{\mathop{\sup}_{t\in[0,T)}\,\mathop{\sup}_{x\in M}\,\mu_t(x)<\infty}$?  
Furthermore, if such a pinching condition were found, 
does the mean curvature flow starting from a curvature-adapted immersion $f$ satisfying $\widehat S=0$ and the pinching condition asymptote to 
the mean curvature flow starting from a curvature-adapted equifocal hypersurface?  Hence, does the flow collapse to a focal submanifold 
of the curvature-adapted equifocal hypersurface?  (see \cite{TT} about the notion of an equifocal hypersurface).
}

\vspace{0.5truecm}

If the above questions are solved affirmatively, then we can derive that there are only finitely many of diffeomorphism classes of 
compact hypersurfaces in a simply connected compact symmetric space satisfying $\widehat S=0$ and the pinching condition 
of the norms $\|A\|,\,\|\nabla A\|$ and $\|\nabla^2A\|$ by using the finiteness theorem for curvature-adapted equifocal hypersurfaces 
in a simply connected compact symmetric space by J. Q. Ge and C. Qian (\cite{GQ}).  

In Section 2, we recall the evolution equations for the basic geometric quantities along the mean curvature flow.  
In Section 3, we derive the evolution equation for the normal Jacobi operator $\widetilde R(\xi_t)$.  
In Section 4, we prove Theorems A, B and D.  

\section{Evolution equations} 
Let $M$ and $(\widetilde M,\widetilde g)$ be as in Introduction.  
Assume that $\widetilde M$ is a locally symmetric space.  
Denote by $\widetilde{\nabla}$ the Levi-Civita connection of $\widetilde g$.  
Also, denote by $\widetilde R, \widetilde{\rm Ric}$ and ${\widetilde R}^S$ the curvature tensor, the Ricci tensor 
and the scalar curvature of $\widetilde g$.  
Let $\{f_t\}_{t\in[0,T)}$ be the mean curvature flow starting from $f$.  
Define a map $F:M\times[0,T)\to\widetilde M$ by 
$F(x,t):=f_t(x)$ ($(x,t)\in M\times[0,T)$).  
Let $g_t,\nabla^t,R_t,\xi_t,A_t$ and $H_t$ be as in Introduction, and $h_t$ be the second fundamental form 
of $f_t$ for $-\xi_t$.  
Also, let $\pi_M$ be the projection of $M\times[0,T)$ onto $M$.  
For a vector bundle $E$ over $M$, denote by $\pi_M^{\ast}E$ the induced bundle 
of $E$ by $\pi_M$.  Also denote by $E_x$ the fibre of $E$ over $x$ and by $\Gamma(E)$ the space of all 
sections of $E$.  Define a section $g$ of $\pi_M^{\ast}(T^{(0,2)}M)$ by 
$g(x,t)=(g_t)_x$ ($(x,t)\in M\times[0,T)$), where $T^{(0,2)}M$ is the 
$(0,2)$-tensor bundle of $M$.  Similarly, we define a section $R$ of $\pi_M^{\ast}(T^{(1,3)}M)$, 
$h$ of $\pi_M^{\ast}(T^{(0,2)}M)$, a section $A$ of $\pi_M^{\ast}(T^{(1,1)}M)$, 
a map $H:M\times[0,T)\to{\Bbb R}$ and a section $\xi$ of $F^{\ast}T\widetilde M$ in terms of 
$h_t,\ A_t,\,\,H_t$ and $\xi_t$, respectively.  
The bundle $\pi_M^{\ast}(TM)$ is regarded as a subbundle of $T(M\times[0,T))$ under the identification of 
$((x,t),v)\in(\pi_M^{\ast}TM)_{(x,t)}=\{(x,t)\}\times T_xM$ and $v^L_{(x,t)}(\in T_{(x,t)}(M\times[0,T))$, 
where $v^L_{(x,t)}$ is the horizontal lift of $v$ to $(x,t)$ with respect to $\pi_M$.  
Also, the fibre $\pi_M^{\ast}(TM)_{(x,t)}$ is identified with $T_xM$ under the identification of 
$((x,t),v)(\in\pi_M^{\ast}(TM)_{(x,t)})$ and $v$.  
For a section $B$ of $\pi_M^{\ast}(T^{(r,s)}M)$, 
we define $\displaystyle{\frac{\partial B}{\partial t}}$ by 
$\displaystyle{\left(\frac{\partial B}{\partial t}\right)_{(x,t)}
:=\frac{d B_{(x,t)}}{dt}}$, where the right-hand side of this relation is 
the derivative of the vector-valued function $t\mapsto B_{(x,t)}\,(\in 
T^{(r,s)}_xM)$.  
For a tangent vector field $X$ on $M$ (or an open set $U$ of $M$), we define 
a section $\overline X$ of $\pi_M^{\ast}TM$ (or $\pi_M^{\ast}TM\vert_U$) by $\overline X_{(x,t)}:=((x,t),X_x)$ 
($(x,t)\in M\times[0,T)\,\,({\rm or}\,\,U\times[0,T))$).  
Define a connection $\nabla$ of $\pi_M^{\ast}TM$ by 
$$(\nabla_vX)_{(\cdot,t)}:=\nabla^t_vX_{(\cdot,t)}\,\,\,{\rm and}\,\,\,
\nabla_{\frac{\partial}{\partial t}}X:=\frac{dX_{(x,\cdot)}}{dt}$$
for $v\in(\pi_M^{\ast}TM)_{(x,t)}(=T_xM)$ and $X\in\Gamma(\pi_M^{\ast}TM)$, where 
$\displaystyle{\frac{dX_{(x,t)}}{dt}}$ is the derivative of the 
vector-valued function $t\mapsto X_{(x,t)}\,(\in T_xM)$.  
Let $\{{\mathcal S}_t\}_{t\in[0,T)}$ be a $C^{\infty}$-family of a $(r,s)$-tensor fields on $M$ 
and ${\mathcal S}$ a section of $\pi_M^{\ast}(T^{(r,s)}M)$ defined by 
${\mathcal S}_{(x,t)}:=({\mathcal S}_t)_x$.  We define a section $\triangle\,{\mathcal S}$ of $\pi_M^{\ast}(T^{(r,s)}M)$ by 
$$(\triangle\,{\mathcal S})_{(x,t)}:=\sum_{i=1}^n\nabla_{e_i}\nabla_{e_i}{\mathcal S},$$
where $\nabla$ is the connection of $\pi_M^{\ast}(T^{(r,s)}M)$ (or 
$\pi_M^{\ast}(T^{(r,s+1)}M)$) induced from $\nabla$ and $\{e_1,\cdots,e_n\}$ 
is an orthonormal base of $T_xM$ with respect to $(g_t)_x$.  
Also, we define a section $\triangle_t{\mathcal S}_t$ of $T^{(r,s)}M$ by 
$$(\triangle_t{\mathcal S}_t)_x:=(\triangle {\mathcal S})_{(x,t)}\quad\,\,(x\in M).$$
Let $E$ be a vector bundle over $M$.  For a section ${\mathcal S}$ of 
$\pi_M^{\ast}(T^{(0,r)}M\otimes E)$, we define 
$\displaystyle{{\rm Tr}_g^{\bullet}\,{\mathcal S}(\cdots,\mathop{\bullet}^j,
\cdots,\mathop{\bullet}^k,\cdots)}$ by 
$$({\rm Tr}_g^{\bullet}\,{\mathcal S}(\cdots,\mathop{\bullet}^j,\cdots,\mathop{\bullet}^k,\cdots))_{(x,t)}
=\sum_{i=1}^n{\mathcal S}_{(x,t)}(\cdots,\mathop{e_i}^j,\cdots,\mathop{e_i}^k,\cdots)$$
$((x,t)\in M\times[0,T))$, where $\{e_1,\cdots,e_n\}$ is an orthonormal base of $T_xM$ with respect to $(g_t)_x$, 
$\displaystyle{{\mathcal S}(\cdots,\mathop{\bullet}^j,\cdots,\mathop{\bullet}^k,\cdots)}$ means that 
$\bullet$ is entried into the $j$-th component and 
the $k$-th component of ${\mathcal S}$ and $\displaystyle{{\mathcal S}_{(x,t)}(\cdots,\mathop{e_i}^j,\cdots,\mathop{e_i}^k,\cdots)}$ 
means that $e_i$ is entried into the $j$-th component and the $k$-th component of ${\mathcal S}_{(x,t)}$.  
By using the normal Jacobi operator $\widetilde R(\xi_t)$ of $f_t$ we define a section $\widetilde R(\xi)$ of 
$\pi_M^{\ast}T^{(1,1)}M$ by $\widetilde R(\xi)_{(x,t)}:=\widetilde R(\xi_t)_x(=((x,t),\widetilde R(\xi_t)_x))\,\,\,((x,t)\in M\times[0,T))$.  
Also, by using the $(1,1)$-tensor field $\widetilde{\rm Ric}^{\sharp}$ over $\widetilde M$ 
(which is define by $\widetilde g(\widetilde{\rm Ric}^{\sharp}(X),Y)=\widetilde{\rm Ric}(X,Y)\,\,\,(X,Y\in T\widetilde M)$), 
we define a $(1,1)$-tensor field $(\widetilde{\rm Ric}^{\sharp})^T_t$ over $M$ by 
$(\widetilde{\rm Ric}^{\sharp})^T_t=f_{t\ast}^{-1}\circ{\rm pr}^T_t\circ\widetilde{\rm Ric}^{\sharp}\circ f_{t\ast}$, where 
${\rm pr}^T_t$ is the orthogonal projection of $f_t^{\ast}T\widetilde M$ onto $f_{t\ast}(TM)$.  
Denote by $(\widetilde{\rm Ric}^{\sharp})^T$ the section of $\pi_M^{\ast}T^{(1,1)}M$ defined by using 
$(\widetilde{\rm Ric}^{\sharp})^T_t$.  
Similarly, we define a $(1,3)$-tensor field ${\widetilde R}^T_t$ over $M$ by 
${\widetilde R}^T_t=f_{t\ast}^{-1}\circ{\rm pr}^T_t\circ \widetilde R\circ(f_{t\ast}\times f_{t\ast}\times f_{t\ast})$.  
Denote by ${\widetilde R}^T$ the section of $\pi_M^{\ast}T^{(1,3)}M$ defined by using ${\widetilde R}^T_t$.  
Since $\widetilde M$ is locally symmetric and irreducible, 
$\widetilde{\nabla}\widetilde R=0$ holds and it is Einstein, that is, 
$$\widetilde{\rm Ric}^{\sharp}=\frac{{\widetilde R}^S}{n}{\rm id}.\leqno{(2.1)}$$

According to (i) of Lemma 3.3 in [Hu2], we have the following evolution equation.  

\vspace{0.5truecm}

\noindent
{\bf Lemma 2.1.} {\sl The family $\{g_t\}_{t\in[0,T)}$ satisfies the following evolution equation:
$$\frac{\partial g}{\partial t}=-2H_th_t.$$
}

\vspace{0.5truecm}

According to (ii) of Lemma 3.3 in [Hu2], we have the following evolution equation.  

\vspace{0.5truecm}

\noindent
{\bf Lemma 2.2.} {\sl The family $\{\xi_t\}_{t\in[0,T)}$ satisfies the following evolution equation:
$$\frac{\partial\xi}{\partial t}=-F_{\ast}({\rm grad}_{g_t}H_t),$$
where ${\rm grad}_{g_t}H_t$ is the element of $\pi_M^{\ast}(TM)$ such that 
$dH_t(X)=g_t({\rm grad}_{g_t}H_t,X)$ for any $X\in\pi_M^{\ast}(TM)$.}

\vspace{0.5truecm}

According to (i) of Lemma 3.3 and Theorem 3.4 in [Hu2], we have the following evolution equation.  

\vspace{0.5truecm}

\noindent
{\bf Lemma 2.3.} {\sl The family $\{A_t\}_{t\in[0,T)}$ satisfies the following evolution equation:
$$\begin{array}{l}
\displaystyle{\frac{\partial A}{\partial t}
=\triangle_tA_t+\left({\rm Tr}(A_t^2)+{\rm Tr}\,\widetilde R(\xi_t)\right)A_t
+2A_t^3-2{\rm Tr}(A_t^2)A_t-\frac{2{\widetilde R}^S}{n}A_t}\\
\hspace{1truecm}\displaystyle{+A_t\circ \widetilde R(\xi_t)+\widetilde R(\xi_t)\circ A_t
+2{\rm Tr}^{\bullet}_{g_t}R_t(\cdot,\bullet)(A_t(\bullet)).}
\end{array}$$
}

\vspace{0.5truecm}

\noindent
{\it Proof.} According to (i) of Lemma 3.3 and Theorem 3.4 in [Hu2] and $(2.1)$, we have 
$$\begin{array}{l}
\displaystyle{\frac{\partial A}{\partial t}
=\triangle_tA_t+\left({\rm Tr}(A_t^2)+{\rm Tr}\,\widetilde R(\xi_t)\right)A_t
-\frac{2{\widetilde R}^S}{n}A_t}\\
\hspace{1truecm}\displaystyle{+A_t\circ \widetilde R(\xi_t)+\widetilde R(\xi_t)\circ A_t
+2{\rm Tr}^{\bullet}_{g_t}{\widetilde R}^T(\cdot,\bullet)(A_t(\bullet)).}
\end{array}$$
On the other hand, according to the Gauss equation, we have 
$${\widetilde R}^T(X,Y)Z=R_t(X,Y)Z+h_t(X,Z)A_tY-h_t(Y,Z)A_tX\leqno{(2.2)}$$
for any tangent vector fields $X,Y$ and $Z$ on $M$.  
From these relations, we can derive the desired evolution equation.  \qed
%\begin{flushright}q.e.d.\end{flushright}

%\noindent
%{\bf Lemma 2.3.} {\sl The family $\{h_t\}_{t\in[0,T)}$ satisfies the following evolution equation:
%$$\begin{array}{l}
%\displaystyle{\frac{\partial h}{\partial t}(X,Y)
%=(\triangle h)(X,Y)-2Hh(AX,Y)+{\rm Tr}\left({\rm Tr}(A^2)+{\rm Tr}\,R(\xi)\right)h(X,Y)}\\
%\hspace{2.25truecm}\displaystyle{-\widetilde{\rm Ric}(X,AY)+h(R(\xi)(X),Y)-\widetilde{\rm Ric}(Y,AX)+h(R(\xi)(Y),X)}\\
%\hspace{2.25truecm}\displaystyle{+2{\rm Tr}^{\bullet}_{g}R(A\bullet,X,Y,\bullet)}\\
%\hspace{8.9truecm}\displaystyle{(X,Y\in\pi_M^{\ast}TM).}
%\end{array}$$
%}

\vspace{0.5truecm}

According to Corollary 3.5 in [Hu2], we have the following evolution equation.  

\vspace{0.5truecm}

\noindent
{\bf Lemma 2.4.} {\sl The family $\{H_t\}_{t\in[0,T)}$ satisfies the following evolution equation:
$$\begin{array}{l}
\displaystyle{\frac{\partial H}{\partial t}
=\triangle_t\,H_t+\left({\rm Tr}(A_t^2)+{\rm Tr}\,\widetilde R(\xi_t)\right)H_t.}
\end{array}$$
}

\section{Evolution of the normal Jacobi operator} 
We use the notations in Sections 1 and 2.  
Assume that $\widetilde M$ is a locally symmetric space.  
In this section, we derive the evolution equation for the family $\{\widetilde R(\xi_t)\}_{t\in[0,T)}$ 
of the normal Jacobi operators.  
%By using the $(1,2)$-tensor fields $(\widetilde R_1(\xi_t))_t^T$ and $(\widetilde R_3(\xi_t))_t^T$ over $M$ by 
%$$\begin{array}{l}
%\displaystyle{\widetilde R_1(\xi_t)_t^T:=f_{t\ast}^{-1}\circ{\rm pr}^T_t\circ \widetilde R(\xi_t,\cdot)\cdot\circ
%(f_{t\ast}\times f_{t\ast}),}\\
%\displaystyle{\widetilde R_3(\xi_t)_t^T:=f_{t\ast}^{-1}\circ{\rm pr}^T_t\circ \widetilde R(\cdot,\cdot)\xi_t\circ
%(f_{t\ast}\times f_{t\ast}),}
%\end{array}$$respctively.  
Denote by $S$ the section of $\pi_M^{\ast}T^{(1,1)}M$ defined by using $S_t$'s and 
by $\widetilde R_i(\xi)$ ($i=1,3$) the sections of $\pi_M^{\ast}T^{(1,2)}M$ defined by using 
$\widetilde R_i(\xi_t)$'s.  
Denote by ${\widetilde{\nabla}}^F$ (resp. ${\widetilde{\nabla}}^{f_t}$) the pull-back connection of 
$\widetilde{\nabla}$ by $F$ (resp. $f_t$).  
First we prepare the following lemma.  

\vspace{0.5truecm}

\noindent
{\bf Lemma 3.1.} {\sl Let $\{Z_t\}_{t\in[0,T)}$ be a $C^{\infty}$-family of tangent vector fields on $M$.  
Then we have 
$$\frac{\partial f_{t\ast}(Z_t)}{\partial t}=f_{t\ast}\left(\frac{\partial Z_t}{\partial t}\right)
-(Z_tH_t)\xi_t-H_tf_{t\ast}(A_tZ_t).$$
}

\vspace{0.5truecm}

\noindent
{\it Proof.} Fix $x_0\in M$.  
Let $\{\phi^t_s\}_{s\in I}$ be the local one-parameter transformation group of $Z_t$ and 
define a map $\widetilde{\delta}:[0,T)^2\times I\to\widetilde M$ by 
$\widetilde{\delta}(t,u,s):=f_t(\phi^u_s(x_0))$ and $\delta:[0,T)\times I\to\widetilde M$ by 
$\delta(t,s):=\widetilde{\delta}(t,t,s)$.  
Then we have 
$$
%\begin{array}{l}\displaystyle{
\frac{\partial\delta}{\partial t}
=\left.\left(\frac{\partial\widetilde{\delta}}{\partial t}
+\frac{\partial\widetilde{\delta}}{\partial u}\right)\right\vert_{u=t}
%}\\\hspace{0.8truecm}\displaystyle{
=-(H_t)_{\phi^t_s(x_0)}(\xi_t)_{\phi^t_s(x_0)}
+\left.f_{t\ast}\left(\frac{\partial\phi^u_s(x_0)}{\partial u}\right)\right\vert_{u=t}.
%}\end{array}
\leqno{(3.1)}$$
Denote by ${\widetilde{\nabla}}^{\delta}$ the pull-back connection of $\widetilde{\nabla}$ by $\delta$.  
Then we have 
\begin{align*}
\left(\frac{\partial f_{t\ast}(Z_t)}{\partial t}\right)_{(x_0,t_0)}
&=\left({\widetilde{\nabla}}^F_{\frac{\partial}{\partial t}}F_{\ast}Z\right)_{(x_0,t_0)}
=\left.\left({\widetilde{\nabla}}^{\delta}_{\frac{\partial}{\partial t}}
\left.\frac{\partial\delta}{\partial s}\right\vert_{s=0}\right)\right\vert_{t=t_0}\\
&=\left.\left({\widetilde{\nabla}}^{\delta}_{\frac{\partial}{\partial s}}
\left.\frac{\partial\delta}{\partial t}\right\vert_{t=t_0}\right)\right\vert_{s=0}.
\end{align*}
By substituting $(3.1)$ into this relation, we can derive 
$$\begin{array}{l}
\displaystyle{\left(\frac{\partial f_{t\ast}(Z_t)}{\partial t}\right)_{(x_0,t_0)}
={\widetilde{\nabla}}^{f_{t_0}}_{(Z_{t_0})_{x_0}}(-H_{t_0}\xi_{t_0})
+\left.\left({\widetilde{\nabla}}^{\widetilde{\delta}}_{\frac{\partial}{\partial u}}
\frac{\partial\widetilde{\delta}}{\partial s}\right)\right\vert_{t=u=t_0,s=0}}\\
\displaystyle{=-((Z_{t_0})_{x_0}H_{t_0})(\xi_{t_0})_{x_0}
-(H_{t_0})_{x_0}f_{t_0\ast}(A_{t_0}(Z_{t_0})_{x_0})
+f_{t_0\ast}\left(\left.\frac{d(Z_u)_{x_0}}{du}\right\vert_{u=t_0}\right).}
\end{array}$$
Therefore, the desired relation follows from the arbitrariness of $(x_0,t_0)$.  \qed
%\begin{flushright}q.e.d.\end{flushright}

\vspace{0.5truecm}

Since $(\widetilde M,\widetilde g)$ is an Einstein space, we have the following relation.  

\vspace{0.5truecm}

\noindent
{\bf Lemma 3.2.} {\sl The following relation holds:
$${\rm grad}_{g_t}H_t={\rm Tr}_{g_t}^{\bullet}(\nabla^t_{\bullet}A_t)(\bullet).$$
}

\vspace{0.5truecm}

\noindent
{\it Proof.} According to the Codazzi equation, we have 
$$\widetilde R_3(\xi_t)(X,Y)=-(\nabla^t_XA_t)(Y)+(\nabla^t_YA_t)(X)\quad\,\,(X,Y\in TM).\leqno{(3.2)}$$
From this relation and the Einsteinity of $(\widetilde M,\widetilde g)$, we obtain 
\begin{align*}
{\rm grad}_{g_t}H_t&=(\bullet\,\mapsto\,{\rm Tr}(\nabla^t_{\bullet}A_t))^{\sharp}
={\rm Tr}_{g_t}^{\bullet}(\nabla^t_{\bullet}A_t)(\bullet)+\widetilde{\rm Ric}^{\sharp}(\xi_t)_T\\
&={\rm Tr}_{g_t}^{\bullet}(\nabla^t_{\bullet}A_t)(\bullet)+\left(\frac{\widetilde R^S}{n}\,\xi_t\right)_T
={\rm Tr}_{g_t}^{\bullet}(\nabla^t_{\bullet}A_t)(\bullet).
\end{align*}
\qed

\vspace{0.5truecm}

By using these lemmas, we can derive the following evolution.  

\vspace{0.5truecm}

\noindent
{\bf Proposition 3.3.} 
{\sl The family $\{\widetilde R(\xi_t)\}_{t\in[0,T)}$ satisfies the following evolution equation:
$$\begin{array}{l}
\displaystyle{\frac{\partial \widetilde R(\xi)}{\partial t}
=\triangle_t\widetilde R(\xi_t)+H_tS_t-\widetilde R(\xi_t)\circ A_t^2+2{\rm Tr}(A_t^2)\widetilde R(\xi_t)}\\
\hspace{1.55truecm}\displaystyle{-2\widetilde R_3(\xi_t)(\cdot,{\rm Tr}^{\bullet}_{g_t}(\nabla^t_{\bullet}A_t)(\bullet))
%}\\\hspace{1.5truecm}\displaystyle{
+2\widetilde R_1(\xi_t)(\cdot,{\rm Tr}^{\bullet}_{g_t}(\nabla^t_{\bullet}A_t)(\bullet))}\\
\hspace{1.55truecm}\displaystyle{-2{\rm Tr}^{\bullet}_{g_t}R_t(\cdot,A_t(\bullet))A_t(\bullet)
-2{\rm Tr}(A_t^3)A_t+2A_t^4.}
\end{array}$$
}

\vspace{0.5truecm}

\noindent
{\it Proof.} 
Take $X$ be a tangent vector field on $M$ and $\overline X$ be the section $\pi_M^{\ast}TM$ defined by 
$\overline X_{(x,t)}:=((x,t),X_x)$.  By using Lemma 3.1 and $\frac{\partial\overline X}{\partial t}=0$, 
we can show 
$$\begin{array}{l}
\displaystyle{\frac{\partial(f_{t\ast}\circ \widetilde R(\xi_t))}{\partial t}(X)
=\frac{\partial f_{t\ast}(\widetilde R(\xi_t)(\overline X))}{\partial t}}\\
\displaystyle{=f_{t\ast}\left(\frac{\partial \widetilde R(\xi_t)}{\partial t}(X)\right)
-((\widetilde R(\xi_t)(X))H_t)\xi_t-H_tf_{t\ast}(A_t(\widetilde R(\xi_t)(X)))}\\
\displaystyle{\equiv f_{t\ast}\left(\frac{\partial \widetilde R(\xi_t)}{\partial t}(X)\right)
-H_tf_{t\ast}(A_t(\widetilde R(\xi_t)(X)))\qquad\,\,({\rm mod}\,\,{\rm Span}\{\xi_t\}).}
\end{array}\leqno{(3.3)}$$
On the other hand, by using Lemmas 2.2, 3.1, $\frac{\partial\overline X}{\partial t}=0$ and 
$\widetilde{\nabla}\widetilde R=0$, we can show 
$$\begin{array}{l}
\displaystyle{\frac{\partial(f_{t\ast}\circ \widetilde R(\xi_t))}{\partial t}(X)
=\frac{\partial F_{\ast}(\widetilde R(\xi)(\overline X))}{\partial t}
=\frac{\partial}{\partial t}\widetilde R(F_{\ast}(\overline X),\xi)\xi}\\
\displaystyle{=\widetilde R\left(\frac{\partial f_{t\ast}(X)}{\partial t},\xi_t\right)\xi_t
+\widetilde R\left(f_{t\ast}(X),\frac{\partial\xi}{\partial t}\right)\xi_t
+\widetilde R\left(f_{t\ast}(X),\xi_t\right)\frac{\partial\xi}{\partial t}}\\
\displaystyle{=-\widetilde R\left((XH_t)\xi_t+H_tf_{t\ast}(A_tX),\xi_t\right)\xi_t
-\widetilde R\left(f_{t\ast}(X),f_{t\ast}({\rm grad}_{g_t}H_t)\right)\xi_t}\\
\hspace{0.5truecm}\displaystyle{-\widetilde R\left(f_{t\ast}(X),\xi_t\right)f_{t\ast}({\rm grad}_{g_t}H_t)}\\
\displaystyle{=-f_{t\ast}\left(H_t\widetilde R(\xi_t)(A_tX)+\widetilde R_3(\xi_t)(X,{\rm grad}_{g_t}H_t)
-\widetilde R_1(\xi_t)(X,{\rm grad}_{g_t}H_t)\right).}
\end{array}\leqno{(3.4)}$$
From $(3.3),\,(3.4)$ and Lemma 3.2, we derive 
$$\frac{\partial \widetilde R(\xi_t)}{\partial t}(X)
=H_tS_t(X)-\widetilde R_3(\xi_t)(X,{\rm Tr}^{\bullet}_g(\nabla_{\bullet}A)(\bullet))
+\widetilde R_1(\xi_t)(X,{\rm Tr}^{\bullet}_g(\nabla_{\bullet}A)(\bullet)).
\leqno{(3.5)}$$

Fix $(x_0,t_0)\in M\times[0,T)$.  Take any $v,w\in T_{x_0}M$.  Let $\widetilde w$ be a tangent vector field on 
a neighborhood of $x_0$ in $M$ with $\widetilde w_{x_0}=w$ and $(\nabla^{t_0}\widetilde w)_{x_0}=0$.  
Then we have 
$${\widetilde{\nabla}}^{f_{t_0}}_v(f_{t_0\ast}(\widetilde R(\xi_{t_0})(\widetilde w)))
=f_{t_0\ast}((\nabla^{t_0}_v\widetilde R(\xi_{t_0}))(w))-h_{t_0}(v,\widetilde R(\xi_{t_0})(w))(\xi_{t_0})_{x_0}.\leqno{(3.6)}$$
On the other hand, from $\widetilde{\nabla}\widetilde R=0$, we have 
$$\begin{array}{l}
\hspace{0.5truecm}\displaystyle{{\widetilde{\nabla}}^{f_{t_0}}_v(f_{t_0\ast}(\widetilde R(\xi_{t_0})(\widetilde w)))
={\widetilde{\nabla}}^{f_{t_0}}_v(\widetilde R(f_{t_0\ast}(\widetilde w),\xi_{t_0})\xi_{t_0})}\\
\displaystyle{=\widetilde R({\widetilde{\nabla}}^{f_{t_0}}_vf_{t_0\ast}(\widetilde w),(\xi_{t_0})_{x_0})(\xi_{t_0})_{x_0}
+\widetilde R(f_{t_0\ast}(w),{\widetilde{\nabla}}^{f_{t_0}}_v\xi_{t_0})\xi_{t_0}}\\
\hspace{0.5truecm}\displaystyle{+\widetilde R(f_{t_0\ast}(w),\xi_{t_0}){\widetilde{\nabla}}^{f_{t_0}}_v\xi_{t_0}}\\
\displaystyle{=-h(v,w)\widetilde R((\xi_{t_0})_{x_0},(\xi_{t_0})_{x_0})(\xi_{t_0})_{x_0}
+\widetilde R(f_{t_0\ast}(w),f_{t_0\ast}(A_{t_0}(v)))\xi_{t_0}}\\
\hspace{0.5truecm}\displaystyle{+\widetilde R(f_{t_0\ast}(w),\xi_{t_0})f_{t_0\ast}(A_{t_0}(v))}\\
\displaystyle{=\widetilde R(f_{t_0\ast}(w),f_{t_0\ast}(A_{t_0}(v)))\xi_{t_0}+\widetilde R(f_{t_0\ast}(w),
\xi_{t_0})f_{t_0\ast}(A_{t_0}(v)).}
\end{array}\leqno{(3.7)}$$
From $(3.6)$ and $(3.7)$, we can derive 
$$(\nabla^{t_0}_v\widetilde R(\xi_{t_0}))(w)=\widetilde R_3(\xi_{t_0})(w,A_{t_0}(v))-\widetilde R_1(\xi_{t_0})(w,A_{t_0}(v)).$$
From the arbitrarinesses of $v,w$ and $(x_0,t_0)$, we have 
$$(\nabla^t_X\widetilde R(\xi_t))(Y)=\widetilde R_3(\xi_t)(Y,A_tX)-\widetilde R_1(\xi_t)(Y,A_tX)\leqno{(3.8)}$$
for any $X,Y\in\Gamma(TM)$ and any $t\in[0,T)$.  
Let $v,w$ and $\widetilde w$ be as above.  Also, let $\widetilde v$ be a tangent vector field on 
a neighborhood of $x_0$ in $M$ with $\widetilde v_{x_0}=v$ and $(\nabla^{t_0}\widetilde v)_{x_0}=0$.   
Then we have 
$$\begin{array}{l}
\displaystyle{f_{t_0\ast}\left((\nabla_v^{t_0}\nabla_v^{t_0}\widetilde R(\xi_{t_0}))(w)\right)
=f_{t_0\ast}\left(\nabla_v^{t_0}((\nabla_{\widetilde v}^{t_0}\widetilde R(\xi_{t_0}))(\widetilde w))\right)}\\
\displaystyle{\equiv\left({\widetilde{\nabla}}^{f_{t_0}}_v
f_{t_0\ast}((\nabla_{\widetilde v}^{t_0}\widetilde R(\xi_{t_0}))(\widetilde w))\right)_T\qquad\,\,
({\rm mod}\,\,{\rm Span}\{(\xi_{t_0})_{x_0}\}).}
%+h_{t_0}(v,(\nabla_v^{t_0}\widetilde R(\xi_{t_0}))(w))(\xi_{t_0})_{x_0}.}
%\\\displaystyle{={\widetilde{\nabla}}^{f_{t_0}}_vf_{t_0\ast}(\widetilde R_3(\xi_t)(\widetilde w,A_t(\widetilde v)))}\\
%\hspace{0.5truecm}\displaystyle{+{\widetilde{\nabla}}^{f_{t_0}}_vf_{t_0\ast}
%(-\widetilde R_1(\xi_{t_0})(\widetilde w,A_{t_0}(\widetilde v))+h_{t_0}(\widetilde v,\widetilde R(\xi_{t_0})(\widetilde w))
%\xi_{t_0})}\\\hspace{0.5truecm}\displaystyle{+h_{t_0}(v,(\nabla^{t_0}_v\widetilde R(\xi_{t_0}))(\widetilde w))\xi_t)}
\end{array}\leqno{(3.9)}$$
From $(3.8)$, $(3.9)$ and $\widetilde{\nabla}\widetilde R=0$, we can derive 
$$\begin{array}{l}
\displaystyle{(\nabla_v^{t_0}\nabla_v^{t_0}\widetilde R(\xi_{t_0}))(w)
=h_{t_0}(v,w)\widetilde R(\xi_{t_0})(A_{t_0}v)-2h_{t_0}(v,A_{t_0}v)\widetilde R(\xi_{t_0})(w)}\\
\hspace{3.65truecm}\displaystyle{
%+h_{t_0}(v,\widetilde R(\xi_{t_0})(w))A_{t_0}v
+\widetilde R_3(\xi_{t_0})(w,(\nabla^{t_0}_vA_{t_0})(v))
-\widetilde R_1(\xi_{t_0})(w,(\nabla^{t_0}_vA_{t_0})(v))}\\
\hspace{3.65truecm}\displaystyle{+2R_{t_0}(w,A_{t_0}v)(A_{t_0}v)+2h_{t_0}(A_{t_0}v,A_{t_0}v)A_{t_0}w}\\
\hspace{3.65truecm}\displaystyle{-2h_{t_0}(w,A_{t_0}v)A_{t_0}^2v.}
\end{array}$$
%(somewhat long calculation)
From the arbitrarinesses of $v,w$ and $(x_0,t_0)$, we have 
$$\begin{array}{l}
\displaystyle{(\nabla_X^t\nabla_X^t\widetilde R(\xi_t))(Y)
=h_t(X,Y)\widetilde R(\xi_t)(A_tX)-2h_t(X,A_tX)\widetilde R(\xi_t)(Y)}\\
\hspace{3.5truecm}\displaystyle{
%+h_t(X,\widetilde R(\xi_t)(Y))A_tX
+\widetilde R_3(\xi_t)(Y,(\nabla^t_XA_t)(X))
-\widetilde R_1(\xi_t)(Y,(\nabla^t_XA_t)(X))
}\\
\hspace{3.5truecm}\displaystyle{+2R_t(Y,A_tX)(A_tX)
+2h_t(A_tX,A_tX)A_tY}\\
\hspace{3.5truecm}\displaystyle{-2h_t(Y,A_tX)A_t^2X}
\end{array}\leqno{(3.10)}$$
for any $X,Y\in\Gamma(TM)$ and any $t\in[0,T)$.  
Hence we can derive 
$$\begin{array}{l}
\displaystyle{(\triangle_t\widetilde R(\xi_t))(X)
=\widetilde R(\xi_t)(A_t^2X)-2{\rm Tr}(A_t^2)\widetilde R(\xi_t)(X)}\\
\hspace{2.85truecm}\displaystyle{
%+A_t^2(\widetilde R(\xi_t)(X))
+\widetilde R_3(\xi_t)(X,{\rm Tr}_{g_t}^{\bullet}(\nabla^t_{\bullet}A_t)(\bullet))
-\widetilde R_1(\xi_t)(X,{\rm Tr}_{g_t}^{\bullet}(\nabla^t_{\bullet}A_t)(\bullet))}\\
\hspace{2.85truecm}\displaystyle{+2{\rm Tr}_{g_t}^{\bullet}R_t(X,A_t(\bullet))(A_t(\bullet))
+2{\rm Tr}(A_t^3)A_tX-2A_t^4X}
\end{array}\leqno{(3.11)}$$
for any $X\in\Gamma(TM)$ and any $t\in[0,T)$.  
From $(3.5)$ and $(3.11)$, we can derive the desired evolution equation.  \qed
%\begin{flushright}q.e.d.\end{flushright}

\section{Proofs of Theorems A, B and D} 
In this section, we use the notations in Sections 1-3.  
First we shall prove Theorem A stated in Introduction.  
For its purpose, we shall show the following lemma.  

\vspace{0.5truecm}

\noindent
{\bf Lemma 4.1.} {\sl Assume that $\widetilde M$ is a locally symmetric space and that 
$f(:M\hookrightarrow\widetilde M)$ is curvature-adapted.  Also, let $f^r$ be the end-point map of 
the hypersurface $f(M)$ for $r\xi$ (i.e., $f^r(x):=\exp^{\perp}(r\xi_x)\,\,\,(x\in M)$), 
where $\exp^{\perp}$ is the normal exponential map of $f$ and $r$ is a real number close to zero sufficiently.  
Then the following statements (i) and (ii) hold.  

(i)\ \ $f^r$ also is curvature-adapted.  

(ii)\ \ Furthermore, if both the shape operator and the normal Jacobi operator of $f$ have constant eigenvalues, 
then both the shape operator and the normal Jacobi operator of $f^r$ also have constant eigenvalues.  
}

\vspace{0.5truecm}

\noindent
{\it Proof.} Fix $x\in M$.  
%%Denote by $G$ the isometry group of $\widetilde M$ and by $K$ the isotropy group of $G$ at $f(x)$.  Then we have $\widetilde M=G/K$.  
Let $\gamma_{\xi_x}$ be the normal geodesic of $f(M)$ with $\gamma_{\xi_x}'(0)=\xi_x$ and set $\xi^r_x:=\gamma_{\xi_x}'(r)$.  
Then, since $\widetilde M$ is locally symmetric, it is shown that $\xi^r_x$ is a unit normal vector of $f^r(M)$ at $x$.  
Denote by $A$ the shape operator of $f$ for $-\xi$ and by $\widetilde R(\xi)$ the normal Jacobi operator of $f$.  
Also, denote by $A^r$ the shape operator of $f^r$ for $-\xi^r(:x\mapsto\xi^r_x)$ and by $\widetilde R(\xi^r)$ 
the normal Jacobi operator of $f^r$.  Since $f$ is curvatrue-adapted, there exists a base $(e_1,\cdots,e_n)$ of $T_xM$ satisfying 
$A_xe_i=\lambda_ie_i$ and $\widetilde R(\xi_x)(e_i)=\nu_ie_i$, where $\lambda_i$ and $\nu_i$ are real numbers.  
Take a curve $\alpha_i:(-\varepsilon,\varepsilon)\to M$ with $\alpha_i'(0)=e_i$ and 
define a map $\delta_i:(-\varepsilon,\varepsilon)\times[0,r+\varepsilon)\to\widetilde M$ by 
$\delta_i(s,t):=\gamma_{\xi_{\alpha_i(s)}}(t)$ ($(s,t)\in(-\varepsilon,\varepsilon)\times[0,r+\varepsilon)$).  
Define a vector field $Y_i$ along $\gamma_{\xi_x}$ by $Y_i:=\frac{\partial\delta_i}{\partial s}\vert_{s=0}$.  
Since $Y_i$ is the Jacobi field along $\gamma_{\xi_x}$ with $Y_i(0)=f_{\ast}e_i$ and $Y'(0)=f_{\ast}A_xe_i$, it is described as 
$$Y_i(t)=\left(\cos(t\sqrt{\nu_i})-\frac{\lambda_i\sin(t\sqrt{\nu_i})}{\sqrt{\nu_i}}\right)P_{\gamma_{\xi_x}\vert_{[0,t]}}(f_{\ast}e_i),
\leqno{(4.1)}$$
where $P_{\gamma_{\xi_x}\vert_{[0,t]}}$ is the parallel translation along $\gamma_{\xi_x}\vert_{[0,t]}$.  
Note that, in case of $\nu_i\leq 0$, $\cos(t\sqrt{\nu_i})=\cosh(t\sqrt{-\nu_i})$ and 
$$\frac{\sin(t\sqrt{\nu_i})}{\sqrt{\nu_i}}=
\left\{
\begin{array}{cc}
\displaystyle{\frac{\sinh(t\sqrt{-\nu_i})}{\sqrt{-\nu_i}}} & (\nu_i<0)\\
\displaystyle{t} & (\nu_i=0).
\end{array}\right.$$
Here we used a general description of Jacobi fields in a symmetric space (see Section 3 of [TT] or (1.2) of [K1]).  
In general, the description is valid in a locally symmetric space.  
From $(4.1)$, we have 
$$f^r_{\ast}(e_i)=Y_i(r)
=\left(\cos(r\sqrt{\nu_i})-\frac{\lambda_i\sin(r\sqrt{\nu_i})}{\sqrt{\nu_i}}\right)
P_{\gamma_{\xi_x}\vert_{[0,r]}}(f_{\ast}e_i).
\leqno{(4.2)}$$
and 
$$f^r_{\ast}(A^r_xe_i)=Y'_i(r)
=-\left(\sqrt{\nu_i}\sin(r\sqrt{\nu_i})+\lambda_i\cos(r\sqrt{\nu_i})\right)
P_{\gamma_{\xi_x}\vert_{[0,r]}}(f_{\ast}e_i).\leqno{(4.3)}$$
Hence we obtain 
$$A^r_xe_i=-\frac{\sqrt{\nu_i}\tan(r\sqrt{\nu_i})+\lambda_i}
{1-\lambda_i\tan(r\sqrt{\nu_i})/\sqrt{\nu_i}}\,e_i.\leqno{(4.4)}$$

%%Let $\mathfrak g$ (resp. $\mathfrak k$) be the Lie algebra of $G$ (resp. $K$) and 
%%$\mathfrak g=\mathfrak k+\mathfrak p$ be the canonical decomposition associated with the symmetric pair 
%%$(G,K)$.   
%%Denote by $\exp$ the exponential map of $G$.  
%%%%%%and by ${\rm Exp}$ the exponential map of $G/K$ at $f(x)$.  
%%The space $\mathfrak p$ is identified with the tangent space $T_{f(x)}(G/K)$ of $G/K$ at $f(x)$.  
Since $\widetilde M$ is locally symmetric, it is shown that $P_{\gamma_{\xi_x}\vert_{[0,r]}}$ is equal to the differential $\phi_{\ast f(x)}$ 
of a local isometry $\phi$ of a neighborhood of $f(x)$ onto a neighborhood of $f^r(x)$ 
(see the discussion in Page 208 of [He] in the case of a symmetric space).  
%%where we note that $\xi_x\in T_{f(x)}(G/K)=\mathfrak p\subset\mathfrak g$.  
%%Since $\exp(r\xi_x)$ is an isometry (called a transvection) of $G/K$, it follows from $(4.2)$ that 
Since $\phi$ is a local isometry, it follows from $(4.2)$ that 
$$\begin{array}{l}
\displaystyle{\widetilde R(\xi^r_x)(e_i)=(f^r_{\ast})^{-1}(\widetilde R(f^r_{\ast}(e_i),\xi^r_x)\xi^r_x)}\\
\displaystyle{=\left(\cos(r\sqrt{\nu_i})-\frac{\lambda_i\sin(r\sqrt{\nu_i})}{\sqrt{\nu_i}}\right)}\\
\hspace{0.5truecm}\displaystyle{\times(f^r_{\ast})^{-1}
(\widetilde R(P_{\gamma_{\xi_x}\vert_{[0,r]}}(f_{\ast}e_i),P_{\gamma_{\xi_x}\vert_{[0,r]}}({\xi}_x))
P_{\gamma_{\xi_x}\vert_{[0,r]}}({\xi}_x))}\\
\displaystyle{=\left(\cos(r\sqrt{\nu_i})-\frac{\lambda_i\sin(r\sqrt{\nu_i})}{\sqrt{\nu_i}}\right)}\\
\hspace{0.5truecm}\displaystyle{\times (f^r_{\ast})^{-1}(\widetilde R(\phi_{\ast x}(f_{\ast}e_i),\phi_{\ast x}({\xi}_x))
\phi_{\ast x}({\xi}_x))}\\
\displaystyle{=\left(\cos(r\sqrt{\nu_i})-\frac{\lambda_i\sin(r\sqrt{\nu_i})}{\sqrt{\nu_i}}\right)
%}\\\hspace{0.5truecm}\displaystyle{
\times((f^r_{\ast})^{-1}\circ\phi_{\ast x})(\widetilde R(f_{\ast}e_i,\xi_x)\xi_x)}\\
\displaystyle{=\left(\cos(r\sqrt{\nu_i})-\frac{\lambda_i\sin(r\sqrt{\nu_i})}{\sqrt{\nu_i}}\right)
%}\\\hspace{0.5truecm}\displaystyle{
\times((f^r_{\ast})^{-1}\circ P_{\gamma_{\xi_x}\vert_{[0,r]}}\circ f_{\ast})(\nu_ie_i)}\\
\displaystyle{=\nu_ie_i.}
\end{array}\leqno{(4.5)}$$
Thus $A^r_x$ and $\widetilde R(\xi^r_x)$ are simultaneously diagonalized with respect to $(e_1,\cdots,e_n)$, that is, 
they commute to each other.  Hence $f^r$ is curvature-adapted.  Thus the statement (i) has been proved.  
The statement (ii) also follows from $(4.4)$ and $(4.5)$.  \qed

\vspace{0.5truecm}

By using this lemma, we shall prove Theorem A.  

\vspace{0.5truecm}

\noindent
{\it Proof of Theorem A.} Let $\widetilde M$ and $f$ be as in the statement of Theorem A.  
Since $A$ and $\widetilde R(\xi)$ have constant eigenvalues by the assumption, it follows from (ii) of Lemma 4.1 that 
$A^r$ and $\widetilde R(\xi^r)$ have constant eigenvalues for any constant $r$ sufficiently close to $0$.  
Hence $f$ and $f^r$ are of constant mean curvature.  This implies that $f(M)$ and $f^r(M)$ are locally isoparametric.  
Since $f$ and $f^r$ are of constant mean curvature, we see that the mean curvature flow $\{f_t\}_{t\in[0,T)}$ is described as 
$f_t=f^{r_t}$ for some $C^{\infty}$-correspndence $t\mapsto r_t$.  Hence $f_t$ is curvature-adapted and locally isoparametric by Lemma 4.1.  \qed

\vspace{0.5truecm}

Next we shall prove Theorem B stated in Introduction.  
For its purpose, we first derive the following evolution equation for the family 
$\{S_t=[A_t,\widetilde R(\xi_t)]\}_{t\in[0,T)}$ in terms of Lemma 2.3 and Proposition 3.3.  

\vspace{0.5truecm}

\noindent
{\bf Lemma 4.2.} {\sl The family $\{S_t\}_{t\in[0,T)}$ satisfies the following evolution equation:
$$\begin{array}{l}
\displaystyle{\frac{\partial S}{\partial t}-\triangle_tS_t
=\left({\rm Tr}(A_t^2)+{\rm Tr}\,\widetilde R(\xi_t)-\frac{2{\widetilde R}^S}{n}\right)S_t}\\
\hspace{2.45truecm}\displaystyle{+H_t[A_t,S_t]-S_t\circ A_t^2+[A_t,\widetilde R(\xi_t)^2]+2[A_t^3,\widetilde R(\xi_t)]-\widehat S_t,}
%\hspace{2.45truecm}\displaystyle{-2\left[A,{\rm Tr}^{\bullet}_gR(\cdot,A(\bullet))(A(\bullet))\right]
%-2\left[\widetilde R(\xi),{\rm Tr}^{\bullet}_gR(\cdot,\bullet)(A(\bullet))\right]}\\
%\hspace{2.45truecm}\displaystyle{-2\left[A,\,(\widetilde R_3(\xi)-\widetilde R_1(\xi))
%(\cdot,{\rm Tr}^{\bullet}_g(\nabla_{\bullet}A)(\bullet))\right]}\\
%\hspace{2.45truecm}\displaystyle{-2{\rm Tr}^{\bullet}_g\left[\nabla_{\bullet}A,\,
%\nabla_{\bullet}\widetilde R(\xi)\right].}
\end{array}\leqno{(4.6)}$$
where $\widehat S_t$ is as in $(1.4)$.}

\vspace{0.5truecm}

\noindent
{\it Proof.} Clearly we have 
$$\begin{array}{l}
\displaystyle{\frac{\partial S}{\partial t}-\triangle_tS_t
=\left[\frac{\partial A}{\partial t}-\triangle_tA_t,\widetilde R(\xi_t)\right]
+\left[A_t,\frac{\partial \widetilde R(\xi_t)}{\partial t}-\triangle_t\widetilde R(\xi_t)\right]}\\
\hspace{2.47truecm}\displaystyle{
-2{\rm Tr}_{g_t}^{\bullet}\left[\nabla^t_{\bullet}A_t,\nabla^t_{\bullet}\widetilde R(\xi_t)\right].}
\end{array}$$
By substituting the evolution equations in Lemma 2.3 and Proposition 3.3 into this relation, 
we can derive the desired evolution equation for $\{S_t\}_{t\in[0,T)}$.  \qed
%\begin{flushright}q.e.d.\end{flushright}

\vspace{0.5truecm}

By using this lemma, we prove Theorem B.  

\vspace{0.5truecm}

\noindent
{\it Proof of Thoerem B}\ Since $f(=f_0)$ is curvature-adpated, we have $S_0=0$.  Hence, from the evolution equation $(4.6)$, we obtain 
$\displaystyle{\left.\frac{\partial S}{\partial t}\right|_{t=0}=-\widehat S_0\not=0}$.  Therefore we can derive the statement of Theorem B.  
\qed

\vspace{0.5truecm}

Denote by $P(S_t)$ the $(-1)$-multiple of the right-hand side of $(4.6)$.  
Let $\rho$ be a function over $M\times[0,T)$ defined by using $\rho_t$'s.  
From $(4.6)$, we can derive the following evolution equation for $\{\rho_t\}_{t\in[0,T)}$ directly.  

\vspace{0.5truecm}

\noindent
{\bf Lemma 4.3.} {\sl The family $\{\rho_t\}_{t\in[0,T)}$ satisfies the following evolution equation:
$$\frac{\partial\rho}{\partial t}-\triangle_t\rho_t=2{\rm Tr}(P(S_t)\circ S_t)+2{\rm Tr}{\rm Tr}_{g_t}^{\bullet}
(\nabla^t_{\bullet}S_t\circ\nabla^t_{\bullet}S_t).$$
}

\vspace{0.5truecm}

For $(1,1)$-tensor fields $\Phi$ and $\Psi$ over $M$, 
we denote ${\rm Tr}(\Phi^{\ast_t}\circ\Psi)$ by $\langle\Phi,\Psi\rangle_t$ and 
${\rm Tr}(\Phi^{\ast_t}\circ\Phi)$ by $\vert\vert\Phi\vert\vert_t^2$, where $\Phi^{\ast_t}$ is the adjoint 
operator of $\Phi$ with respect to $g_t$.  
Define $\vert\vert \widetilde R\vert\vert(:M\to{\Bbb R})$ by 
$$\vert\vert \widetilde R\vert\vert(x):=\max\{\vert\vert \widetilde R(v_1,v_2)v_3\vert\vert\,\,\vert\,
v_i\in T_x\widetilde M\,\,{\rm s.t.}\,\,\vert\vert v_i\vert\vert=1\,\,\,(i=1,2,3)\}\quad\,\,(x\in M),$$
%%(this value is independent of the choice $x\in M$), 
where $\vert\vert\bullet\vert\vert:=\sqrt{\widetilde g(\bullet,\bullet)}$.  
Note that $\vert\vert \widetilde R\vert\vert$ is constant in the case where $\widetilde M$ is 
a Riemannian homogeneous space.  
Clearly we have the following inequalities.  

\vspace{0.5truecm}

\noindent
{\bf Lemma 4.4.} {\sl 
{\rm (i)} 
${\rm Tr}(\widetilde R(\xi_t)^k)\leq n\vert\vert \widetilde R\vert\vert^k$ ($k\in{\Bbb N}$), 

{\rm (ii)} 
${\rm Tr}_{g_t}^{\bullet}{\rm Tr}(\widetilde R_i(\xi_t)(\bullet)\circ \widetilde R_i(\xi_t)(\bullet))
\leq n^2\vert\vert\widetilde R\vert\vert^2\quad(i=1,3)$.  
}

\vspace{0.5truecm}

By using Lemmas 4.3 and 4.4, we can derive the following estimate of the functions $\rho_t$.  

\vspace{0.5truecm}

\noindent
{\bf Proposition 4.5.} {\sl Assume that $\widetilde M$ is a locally symmetric space and that 
$$\mathop{\sup}_{t\in[0,T)}\,\mathop{\sup}_{x\in M}\,\mu_t(x)<\infty,$$
where $\mu_t$ is as in $(1.5)$.  Fix any $T_0\in[0,T)$.  Then $\rho_t\,\,(0\leq t\leq T_0)$ are estimated from above as follows:
$$\rho_t\leq\left(\mathop{\max}_{x\in M}\,\rho_0(x)\right)\cdot e^{C_1(T_0)t}\qquad\,\,(0\leq t<T_0),$$
where $C_1(T_0)$ is defined by 
$$C_1(T_0):=4(2n+1)\left(\mathop{\max}_{(x,t)\in M\times[0,T_0]}\vert\vert A_t\vert\vert_t(x)\right)^2
+10n\vert\vert\widetilde R\vert\vert+2\mathop{\sup}_{t\in[0,T_0]}\,\mathop{\sup}_{x\in M}\,\mu_t(x).$$
}

\vspace{0.5truecm}

\noindent
{\it Proof.} Since $S_t$ is skew-symmetric, so is also $\nabla^t_XS_t$ for any $X\in\Gamma(TM)$.  
Hence $S_t^2$ and $(\nabla^t_XS_t)^2$ are non-positive opearators.  Therefore we obtain 
$$\begin{tabular}{c}
$\rho_t=-{\rm Tr}\,S_t^2\geq 0\quad\,\,$ and $\quad\,\,
{\rm Tr}{\rm Tr}_g^{\bullet}(\nabla^t_{\bullet}S_t\circ\nabla^t_{\bullet}S_t)\leq 0.$\\
\end{tabular}\leqno{(4.7)}$$
Hence, from Lemma 4.3, we have 
$$\frac{\partial\rho}{\partial t}-\triangle_t\rho_t\leq 2{\rm Tr}(P(S_t)\circ S_t).\leqno{(4.8)}$$
For simplicity, we set $C_A(T_0):=\max_{(x,t)\in M\times[0,T_0]}\vert\vert A_t\vert\vert_t(x)$.  
%and $C_{\widetilde R}:=\vert\vert\widetilde R\vert\vert$.  

Now we shall calculate $P(S_t)$.  
Clearly we have 
$$(R_t(X,Y)A^k_t)(Z)=R_t(X,Y)(A_t^kZ)-A_t^k(R_t(X,Y)Z)\leqno{(4.9)}$$
%and 
%$$R_t(X,Y)A_t^2=R_t(X,Y)A_t\circ A_t+A_t\circ R_t(X,Y)A_t$$
for tangent vector fields $X,Y$ and $Z$ on $M$.  
From $(2.2)$, we have 
$${\rm Tr}_{g_t}^{\bullet}\widetilde R^T(X,\bullet)\bullet={\rm Tr}_{g_t}^{\bullet}R_t(X,\bullet)\bullet+A_t^2X-H_tA_tX.$$
On the other hand, from $(2.1)$, we have 
$${\rm Tr}_{g_t}^{\bullet}\widetilde R^T(X,\bullet)\bullet=\frac{\widetilde R^S}{n}X-\widetilde R(\xi_t)(X).$$
Hence we have 
$${\rm Tr}_{g_t}^{\bullet}R_t(X,\bullet)\bullet=\frac{\widetilde R^S}{n}X-\widetilde R(\xi_t)(X)-A_t^2X+H_tA_tX.
\leqno{(4.10)}$$
By using $(4.9)$ and $(4.10)$, we can show 
$$\begin{array}{l}
\displaystyle{{\rm Tr}_{g_t}^{\bullet}R_t(X,A_t(\bullet))(A_t(\bullet))
={\rm Tr}_{g_t}^{\bullet}R_t(X,\bullet)(A_t^2(\bullet))}\\
\displaystyle{={\rm Tr}_{g_t}^{\bullet}(R_t(X,\bullet)A_t^2)(\bullet)
+{\rm Tr}_g^{\bullet}A_t^2(R_t(X,\bullet)\bullet)}\\
\displaystyle{={\rm Tr}_{g_t}^{\bullet}(R_t(X,\bullet)A_t^2)(\bullet)+\frac{\widetilde R^S}{n}A_t^2X}\\
\hspace{0.5truecm}\displaystyle{-(A_t^2\circ\widetilde R(\xi_t))(X)-A_t^4X+H_tA_t^3X.}
\end{array}$$
Also, we have 
$$R_t(X,\cdot)A_t^2=R_t(X,\cdot)A_t\,\circ A_t+A_t\circ R_t(X,\cdot)A_t.$$
Hence we have 
$$[A_t,{\rm Tr}_{g_t}^{\bullet}R_t(X,A_t(\bullet))(A_t(\bullet))]
=[A_t^2,{\rm Tr}_{g_t}^{\bullet}(R_t(X,\bullet)A_t)(\bullet)]-A_t^2\circ S_t.\leqno{(4.11)}$$
Also, by using $(4.9)$ and $(4.10)$, we show 
$$\begin{array}{l}
\displaystyle{{\rm Tr}_{g_t}^{\bullet}R_t(X,\bullet)(A_t(\bullet))
={\rm Tr}_{g_t}^{\bullet}(R_t(X,\bullet)A_t)(\bullet)+\frac{\widetilde R^S}{n}A_tX}\\
\hspace{3.75truecm}\displaystyle{-(A_t\circ\widetilde R(\xi_t))(X)-A_t^3X+H_tA_t^2X.}
\end{array}$$
Hence we obtain 
$$\begin{array}{l}
\displaystyle{[\widetilde R(\xi_t),{\rm Tr}_{g_t}^{\bullet}R_t(X,\bullet)(A_t(\bullet))]
=[\widetilde R(\xi_t),{\rm Tr}_{g_t}^{\bullet}(R_t(X,\bullet)A_t)(\bullet)]-\frac{\widetilde R^S}{n}S_t}\\
\hspace{5truecm}\displaystyle{+S_t\circ\widetilde R(\xi_t)+[A_t^3,\widetilde R(\xi_t)]
-H_t[A_t^2,\widetilde R(\xi_t)].}
\end{array}\leqno{(4.12)}$$
In the sequel, we omit the subscript ``$t$''.  
By using $(4.11)$ and $(4.12)$, we can derive 
$$\begin{array}{l}
\displaystyle{P(S)=\widehat S-\left({\rm Tr}(A^2)+{\rm Tr}\,\widetilde R(\xi)\right)S-H[A,S]
-2H[A^2,\widetilde R(\xi)]}\\
\hspace{1.31truecm}\displaystyle{+S\circ A^2-2A^2\circ S+2S\circ\widetilde R(\xi)
%}\\\hspace{1.27truecm}\displaystyle{
-[A,\widetilde R(\xi)^2]}
%\\
%\hspace{1.4truecm}\displaystyle{+2\left[A,{\rm Tr}^{\bullet}_g(R(\cdot,\bullet)(A^2))(\bullet)\right]
%+2\left[\widetilde R(\xi),{\rm Tr}^{\bullet}_g(R(\cdot,\bullet)A)(\bullet)\right]}\\
%\hspace{1.4truecm}\displaystyle{+\left[A,\,(\widetilde R_3(\xi)-\widetilde R_1(\xi))
%(\cdot,{\rm grad}_gH+{\rm Tr}^{\bullet}_g(\nabla_{\bullet}A)(\bullet))\right]}\\
%\hspace{1.4truecm}\displaystyle{+2{\rm Tr}^{\bullet}_g\left[\nabla_{\bullet}A,\,
%\nabla_{\bullet}\widetilde R(\xi)\right].}
\end{array}\leqno{(4.13)}$$
and hence 
$$\begin{array}{l}
\displaystyle{{\rm Tr}(P(S)\circ S)=-\langle\widehat S,S\rangle
+\left({\rm Tr}(A^2)+{\rm Tr}\,\widetilde R(\xi)\right)\rho
-H{\rm Tr}([A,S]\circ S)}\\
\hspace{2.62truecm}\displaystyle{-2H{\rm Tr}([A^2,\widetilde R(\xi)]\circ S)
-{\rm Tr}(A^2\circ S^2)}\\
\hspace{2.62truecm}\displaystyle{+2{\rm Tr}(\widetilde R(\xi)\circ S^2)
-{\rm Tr}([A,\widetilde R(\xi)^2]\circ S).}
\end{array}\leqno{(4.14)}$$

Now we shall estimate ${\rm Tr}(P(S)\circ S)$ from above.  
By using Lemma 4.4, we have 
$$
%\begin{array}{l}\displaystyle{
\left({\rm Tr}(A^2)+{\rm Tr}\,\widetilde R(\xi)\right)\rho
%}\\\displaystyle{
\leq(\vert\vert A\vert\vert^2+n\vert\vert\widetilde R\vert\vert)\rho
\leq(C_A(T_0)^2+n\vert\vert\widetilde R\vert\vert)\rho,
%}\end{array}
\leqno{(4.15)}$$
$${\rm Tr}(([A,S]\circ S))={\rm Tr}(A\circ S^2-S\circ A\circ S)=0,\leqno{(4.16)}$$
$$\begin{array}{l}
\displaystyle{-2H{\rm Tr}([A^2,\widetilde R(\xi)]\circ S)=-2H{\rm Tr}((A\circ S+S\circ A)\circ S)}\\
\displaystyle{=-4H{\rm Tr}(A\circ S^2)\leq 4n\vert\vert A\vert\vert^2\rho\leq 4nC_A(T_0)^2\rho,}
\end{array}\leqno{(4.17)}$$
$$-{\rm Tr}(A^2\circ S^2)\leq\vert\vert A\vert\vert^2\rho\leq C_A(T_0)^2\rho,\leqno{(4.18)}$$
$${\rm Tr}(\widetilde R(\xi)\circ S^2)\leq n\vert\vert\widetilde R\vert\vert\rho\leqno{(4.19)}$$
%$$-{\rm Tr}([A^3,\widetilde R(\xi)]\circ S)
%=-2{\rm Tr}(A^2\circ S^2)-{\rm Tr}((A\circ S)^2)\leq 3\vert\vert A\vert\vert^2\rho
%\leq 3C_A(T_0)^2\rho\leqno{(4.19)}$$
and 
$$\begin{array}{l}
\displaystyle{-{\rm Tr}([A,\widetilde R(\xi)^2]\circ S)
=-{\rm Tr}(A\circ \widetilde R(\xi)^2\circ S)+{\rm Tr}(\widetilde R(\xi)^2\circ A\circ S)}\\
\displaystyle{=-{\rm Tr}(S\circ(A\circ \widetilde R(\xi))\circ \widetilde R(\xi)))
+{\rm Tr}(S\circ(\widetilde R(\xi)\circ A)\circ \widetilde R(\xi))}\\
\hspace{0.5truecm}\displaystyle{-{\rm Tr}((S\circ \widetilde R(\xi))\circ(A\circ \widetilde R(\xi)))
+{\rm Tr}((S\circ \widetilde R(\xi))\circ(\widetilde R(\xi)\circ A))}\\
\displaystyle{=-{\rm Tr}(S^2\circ \widetilde R(\xi))+{\rm Tr}(S\circ \widetilde R(\xi)\circ S)}\\
\displaystyle{=-2{\rm Tr}(S^2\circ \widetilde R(\xi))\leq 2n\vert\vert\widetilde R\vert\vert\rho}
\end{array}\leqno{(4.20)}$$
on $M\times[0,T_0]$.  
Also, we have 
$$-\langle\widehat S,S\rangle\leq 2\left(\mathop{\sup}_{t\in[0,T_0]}\,\mathop{\sup}_{x\in M}\,\mu_t(x)\right)
\cdot\rho
\leqno{(4.21)}$$
on $M\times[0,T_0]$.  
%From the Gauss and Weingarten formulus, we have 
%$${\widetilde R}^T(X,Y)Z=R(X,Y)Z-h(Y,Z)AX+h(X,Z)AY\quad\,\,(X,Y,Z\in\pi_M^{\ast}TM).\leqno{(4.10)}$$
From $(4.14)-(4.21)$, we obtain 
$$2{\rm Tr}(P(S)\circ S)\leq C_1(T_0)\rho\,\,\,\,{\rm on}\,\,M\times[0,T_0],$$
where $C_1(T_0)$ is the positive constant as in the statement of Proposition 4.5.  
Hence, from $(4.8)$, we obtain 
$$\frac{\partial\rho}{\partial t}-\triangle\rho\leq C_1(T_0)\rho\,\,\,\,{\rm on}\,\,M\times[0,T_0].$$
%%Set $\overline{\rho}:=\rho+\frac{2C_2(T_0)}{C_1(T_0)}$.  
%%Then we have 
%%$$\frac{\partial\overline{\rho}}{\partial t}-\triangle\overline{\rho}\leq C_1(T_0)\overline{\rho}.$$
Furthermore, set $\widehat{\rho}_t:=e^{-C_1(T_0)t}\rho_t$.  
Then we have 
$$\frac{\partial\widehat{\rho}}{\partial t}-\triangle\widehat{\rho}\leq 0.$$
Hence, by the maximum principle, we obtain $\widehat{\rho}_t\leq\max\,\widehat{\rho}_0$, which is equivalent 
to the inequality in the statement.  \qed

\vspace{0.5truecm}

\noindent
{\it Proof of Theorem D.} Take any $T_0\in[0,T)$.  Since $f$ is curvature-adapted, we have $\rho_0=0$.  Hence, 
it follows from Proposition 4.5 that $\rho_t=0$ holds for all $t\in [0,T_0]$.  
Therefore, from the arbitrariness of $T_0$, we can conclude $\rho_t=0$ holds for all $t\in [0,T)$.  \qed

\vspace{0.5truecm}

\noindent
{\it Proof of Corollary E.} Since 
$$\mathop{\inf}_{t\in[0,T)}\,\mathop{\min}_{x\in M}\,\langle(\widehat S_t)_x,(S_t)_x\rangle\geq 0$$
by the assumption, we have 
$$\mathop{\sup}_{t\in[0,T)}\,\mathop{\sup}_{x\in M}\,\mu_t(x)\leq 0.$$
Hence we can derive the statement of Corollary E from Theorem D.  \qed

%%preservability theorem of the curvature-adaptedness along the mean curvature flow.  

%\vspace{10pt}
%\begin{figure}[h]
%\centerline{\input{camf-f1.tex}\hspace{5.15truecm}}
%\caption{On the uniformly boundedness of $\{\mu_t\}_{t\in[0,T)}$}
%%\label{figure1}
%\end{figure}
%\vspace{10pt}

%%\vspace{1truecm}
%%\centerline{{\bf References}}
%%\vspace{0.5truecm}

%\vspace{0.5truecm}
%{\small 
%\rightline{Department of Mathematics, Faculty of Science}
%%\rightline{Tokyo University of Science, 1-3 Kagurazaka}
%\rightline{Shinjuku-ku, Tokyo 162-8601 Japan}
%\rightline{(koike@rs.kagu.tus.ac.jp)}
%}

\end{document}